\documentclass[hyperref={colorlinks,citecolor=red,urlcolor=black,bookmarks=false,hypertexnames=true}]{article}

\usepackage{arxiv}

\usepackage[utf8]{inputenc} 
\usepackage[T1]{fontenc}    
\usepackage{url}            
\usepackage{booktabs}       
\usepackage{amsfonts}       
\usepackage{nicefrac}       
\usepackage{microtype}      
\usepackage{lipsum}         
\usepackage{graphicx}
\usepackage{subcaption}
\usepackage{doi}

\usepackage{colortbl}
\usepackage{listings}
\usepackage{bm}
\usepackage{physics}

\usepackage{amsmath,accents}
\usepackage{cleveref}       

\usepackage[ruled,vlined]{algorithm2e}

\usepackage{appendix}
\usepackage{cite}

\usepackage{multirow}
\usepackage{array}

\usepackage{chngcntr}
\usepackage{epigraph}

\SetKwComment{Comment}{/* }{ */}

\newcommand{\ttl}{%
		About rectified sigmoid function for\\
		enhancing the accuracy of\\
		Physics-Informed Neural Networks
}

\title{\ttl}

\date{\today}

\author{ \href{https://orcid.org/0000-0002-4930-1846}{\includegraphics[scale=0.06]{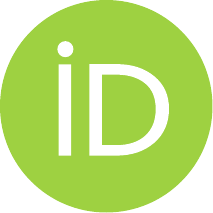}\hspace{1mm}Vasiliy A. Es'kin}\thanks{Corresponding author: Vasiliy Alekseevich Es’kin (\href{vasiliy.eskin@gmail.com}{vasiliy.eskin@gmail.com})} \\
	Department of Radiophysics, University of Nizhny Novgorod\\
	 Nizhny Novgorod, Russia, 603950\\
	 and\\
	 Huawei Nizhny Novgorod Research Center\\
	 Nizhny Novgorod, Russia\\
	\href{vasiliy.eskin@gmail.com}{\texttt{vasiliy.eskin@gmail.com}} \\
	\And
	\href{https://orcid.org/0000-0001-5210-8281}{\includegraphics[scale=0.06]{orcid.pdf}\hspace{1mm}Alexey O. Malkhanov} \\
	Huawei Nizhny Novgorod Research Center\\
	Nizhny Novgorod, Russia\\	\texttt{alexey.malkhanov@gmail.com} \\
	\And
	\href{https://orcid.org/0000-0002-0454-5249}{\includegraphics[scale=0.06]{orcid.pdf}\hspace{1mm}Mikhail E. Smorkalov} \\
	Skolkovo Institute of Science and Technology\\
	Moscow, Russia\\
	and\\
	Huawei Nizhny Novgorod Research Center\\
	Nizhny Novgorod, Russia\\	\texttt{smorkalovme@gmail.com} \\
}

\renewcommand{\headeright}{}
\renewcommand{\undertitle}{}
\renewcommand{\shorttitle}{A preprint}

\renewcommand{\vec}{\bf}

\newcommand{\rsigm}{\text{Re-}\sigma}

\usepackage{tikz}
\newcommand\tikznode[3][]%
{\tikz[remember picture,baseline=(#2.base)]
	\node[minimum size=0pt,inner sep=0pt,#1](#2){#3};%
}

\hypersetup{
pdftitle={A template for the arxiv style},
pdfsubject={math.na, cs.AI, physics.comp-ph},
pdfauthor={Vasiliy A. Es'kin, Danil V. Davydov, Ekaterina D. Egorova, Alexey O. Malkhanov, Mikhail A. Akhukov, Mikhail E. Smorkalov},
pdfkeywords={Deep Learning, Physics-informed Neural Networks, Partial differential equations, Predictive modeling, Computational physics, Nonlinear dynamics},
}

\begin{document}
\maketitle

%

\begin{abstract}
	
	The article is devoted to the study of neural networks with one hidden layer and a modified activation function for solving physical problems. A rectified sigmoid activation function has been proposed to solve physical problems described by the ODE with neural networks. Algorithms for physics-informed data-driven initialization of a neural network and a neuron-by-neuron gradient-free fitting method have been presented for the neural network with this activation function. Numerical experiments demonstrate the superiority of neural networks with a rectified sigmoid function over neural networks with a sigmoid function in the accuracy of solving physical problems (harmonic oscillator, relativistic slingshot, and Lorentz system).

\end{abstract}


\usetikzlibrary {arrows.meta,bending,positioning}

\section{Introduction}

The field of machine learning related to science and engineering~\cite{Ronneberger2021,bodnar2024foundationmodelearth,Raissi2019,Wang2022,Wang2022_2,Lu_2021,Krinitskiy2022,Fanaskov2022,Raissi2017,eskin2023optimal} has seen significant advancements in recent decades, particularly with the integration of physical laws and principles into machine learning models, leading to the development of physics-informed machine learning techniques~\cite{Raissi2019, Wang2022, Wang2022_2}.  In the PINN (Physics-Informed Neural Network) approach, a neural network is trained to approximate the dependences of physical values on spatial and temporal variables for a given physical problem, described by a set of physical equations, together with additional constraints such as initial and boundary conditions. In some cases, data from numerical simulations or experiments can also be used to help train the network. PINNs have been used to solve a wide range of problems described by differential equations in various fields, including thermodynamics, hydrodynamics, mechanics, finance, systems biology, etc.~\cite{Raissi2017,Raissi2019,Kharazmi2021,Yang2021,Cuomo2022,Patel2022,Cai2021,Moseley2020,Yuan2022,daryakenari2023aiaristotle}.

Most current works on the subject using the PINN approach for scientific problems employ small neural networks, with only a few layers, or specific architectures with many layers, which can mitigate the vanishing gradients problem~\cite{wang2024,jiang2024}. However, this issue reduces the representativeness and effectiveness of deep PINNs. In contrast, neural networks with one or two hidden layers avoid this drawback and, despite their simplicity, their potential has not yet been fully explored.

The past article~\cite{eskin2024twohidden} discusses the development of various methods and techniques for initialization and training neural networks with a single hidden layer, as well as the training of a separable physics-informed neural network consisting of neural networks with a single hidden layer to solve physical problems described by ordinary differential equations (ODEs) and partial differential equations (PDEs). Neural networks initialized and trained by the methods proposed in this article demonstrate outstanding accuracy of solutions (at a relatively short learning time) and generalizing properties for parameterized differential equations. In this article~\cite{eskin2024twohidden}, the sigmoid function was taken as the activation function for neural networks. At the same time, we do not think that this activation function is the best for neural networks used to solve physical problems described by differential equations.

This work is devoted to studying neural networks with one hidden layer and modified activation function for solving physical problems. We will briefly list the contributions made in the paper:
\begin{enumerate}
	\item A rectified sigmoid function is proposed for solving physical problems described by the ODE;
	\item Algorithms for physics-informed data-driven (PIDD) initialization of a neural network and gradient-free neuron-by-neuron (NbN) fitting method for a neural network with one hidden layer has been modified with rectified sigmoid function as activation function;	
	\item The accuracy of the ODE solutions of neural networks with a sigmoid activation function is compared with the accuracy of the ODE solutions provided by neural networks with a rectified sigmoid activation function.
\end{enumerate}

The paper is structured as follows. In Section 2, a rectified sigmoid function is proposed for solving physical problems described by the ODE. The algorithms of PIDD initialization of a neural network and NbN fitting method for a neural network with one hidden layer and rectified sigmoid function as activation function are presented. Section 3 shows the results of numerical experiments. Finally, in Section 4 concluding remarks are given.

\section{State of the problem}

We consider a neural network with a single hidden layer for solving problems which are described with ordinary differential equations (ODEs). Consider a system of nonlinear ordinary differential equations, which depend on coordinates $x$, and in general take the following form
\begin{equation}\label{eq1}
	\frac{\partial \vec u}{\partial x} + \mathcal{N} [{\vec u}, x] = 0,\quad x\in [0, X]
\end{equation}
under the initial conditions
\begin{equation}
	{\vec u}(0) = {\vec g}, \label{eq2}
\end{equation}
where ${\vec u}(x)$ denotes the latent solution that is governed by the ODE system of equations (\ref{eq1}) and consists of $n$ components ${\vec u} = (u_1, u_2,\dots,u_n)$, $\mathcal{N}$ is a nonlinear differential operator, ${\vec {g}}$ is initial distribution of ${\vec u}$.

According to PINN approach~\cite{Raissi2019} we approximate the unknown solution ${\vec u}(x)$ with neural networks ${\vec u}_{{\bm \theta}}(x)$, every of which the component ${u}_{{\bm \theta};l}(x)$ is the separated neural network with a single hidden layer, as follows:
\begin{equation}\label{eq3}
	{u}_{{\bm \theta};l}(x) = \sum\limits_{k = 0}^{N-1} {W}^{(2)}_k\sigma\left({W}^{(1)}_k x + {b}^{(1)}_k\right) + {b}^{(2)}_0,
\end{equation}
 where ${W}^{(j)}_k$ is $k$th weight of the $j$th layer, $b^{(j)}_k$ is $k$th bias of the $j$th layer, ${\bm \theta}$ denote all trainable parameters ({weights} and biases) of the neural network ${u}_{{\bm \theta};l}$, $\sigma$ is an activation function, $N$ is number neurons in the hidden layer.
 
 In the article~\cite{eskin2024twohidden}, the sigmoid function was used as an activation function, because it is the most suitable among all the presented activation functions~\cite{kunc2024decadesactivationscomprehensivesurvey} for Euler's method of integrating differential equations. In this paper, we propose using a new activation function, that we call ``rectified sigmoid'', which corresponds better to Euler's approach for integrating differential equations than the original sigmoid function. Rectified sigmoid function is defined as follows
 \begin{equation}\label{eq4}
 	\rsigm(x) = \left\{\begin{array}{ll}
 		0, & x < -1,\\
 		\dfrac{1}{2}(x + 1), & -1 \leq x \leq 1,\\
 		1, & x > 1,
 	\end{array}\right.
 \end{equation}
 This function can be written as a combination of two rectified linear unit (ReLU~\cite{kunc2024decadesactivationscomprehensivesurvey}) functions as follows
\begin{equation}\label{eq5}
	\rsigm(x) = \frac{1}{2}\left[\text{ReLU}(x + 1) - \text{ReLU}(x - 1)\right].
\end{equation}
Perhaps this circumstance is the reason why neural networks with the ReLU activation function demonstrate relatively good results in solving many practical problems~\cite{eskin2024SPINN}.

To verify the accuracy of a neural network with an activation function $\rsigm$, we use the physics-informed data-driven (PIDD) initialization of neural networks algorithm proposed in paper~\cite{eskin2024twohidden}, which is provided in the Appendix~\ref{appA}  for convenience as algorithm~\ref{alg1}. We used this algorithm for a number of physical problems described by equations~(\ref{eq1}). Algorithm~\ref{alg1} was proposed for the sigmoid activation function. Some changes should be made to it when using the activation function $\rsigm$.
Based on the behaviour of the activation function $\rsigm$, we must assume ${\Delta\zeta}$ to be equal to 1. With this choice of ${\Delta\zeta}$, the derivatives of the outputs ($\rsigm'(x)$) of each neuron in the hidden layer do not contribute to those of neighbouring neurons in the same layer. Therefore, the values of $\kappa_k$ in algorithm~\ref{alg1} are equal to 1.

Let us give an equation for calculating the weights ${W}^{(2)}_{m}$. The weights of the hidden layer have already been initialized using algorithm~\ref{alg1} based on the requirements given above.
Find the derivative of the component of solution $u_{\bm \theta;l}(x)$  at the point with the coordinate $x_m$ ($x_m = m \Delta x$, where $\Delta x = X / N$):
\begin{equation}\label{eq6}
	\left.\frac{\partial u_{\bm \theta;l}(x)}{\partial x}\right\vert_{x=x_m} = \sum\limits_{k = 0}^{N-1} \frac{2\Delta \zeta}{\Delta x}{W}^{(2)}_k\rsigm '\left(\frac{2\Delta \zeta}{\Delta x} \left[x_m -x_k\right]\right).
\end{equation}
Here, $\rsigm'(x)$ means a derivate on the argument of the function $\rsigm(x)$. Pay attention $\rsigm'(0) = 1/2$. We have following equation of derivative of $u_{\bm \theta;l}(x)$
\begin{equation}\label{eq7}
	\left.\frac{\partial u_{\bm \theta;l}(x)}{\partial x}\right\vert_{x=x_m} = \frac{\Delta \zeta}{\Delta x}{W}^{(2)}_{m}.
\end{equation}
On the other hand, derivative of $u_{\bm \theta;l}(x)$ must satisfy the Equations (\ref{eq1}). Using (\ref{eq1}) and (\ref{eq7} we have the following equation for the weights ${W}^{(2)}_{m}$
\begin{equation}\label{eq8}
	{W}^{(2)}_{m} =  - \Delta x \mathcal{N}_l [{\vec u}_m, x_{m}] / \Delta \zeta.
\end{equation}
Here values ${\vec u}_m$ is solution ${\vec u}$ in the point $x_m$, $\mathcal{N}_l [{\vec u}, x_{m}]$ is $l$ component of the result of operator action $\mathcal{N} [{\vec u}, x_{m}]$.

As a result, we have the following physics-informed data-driven (PIDD) initialization of neural networks with activation function $\rsigm$ (Algorithm~\ref{alg2}).
\begin{algorithm}[h!]
	\caption{Physics-informed data-driven initialization of neural networks with activation function $\rsigm$}\label{alg2}
	\KwData{$\left\{{\vec u}_k\right\}_{k=0}^{N-1}$ for $\left\{x_k\right\}_{k=0}^{N-1}$ (uniform grid with a step $\Delta x \gets X / N$)}
	\KwResult{Initialized neural network ${u}_{{\bm \theta};l}(x)$ of PINN ${\vec u}_{\bm \theta}$, which consists of $N$ neurons on hidden layer}
	${\Delta\zeta} \gets 1$\;
	\For{$k=0,\dots, N-1$}{
		${W}^{(1)}_k \gets 2 \Delta \zeta / \Delta x$\;
		${b}^{(1)}_k \gets -2 k \Delta \zeta$\;
		${W}^{(2)}_k \gets - \dfrac{\Delta x}{\Delta \zeta} {\mathcal{N}[{\vec u}_{k}, x_{k}]}$\;}
	${b}^{(2)}_0 \gets {u}_{l;0} - \sum\limits_{k = 0}^{N-1} {W}^{(2)}_k\rsigm \left({W}^{(1)}_k x_0 + {b}^{(1)}_k\right)$\;
	where  ${u}_{l;0} = u_l(0)$.
\end{algorithm}

Similarly way, we modify the physical-informed  neuron-by-neuron training for neural networks with a sigmoid activation function (Algorithm~\ref{alg5} in Appendix~\ref{appB}; see the article for details~\cite{eskin2024twohidden}). Neuron-by-neuron training for neural networks with a rectified sigmoid activation is presented in Algorithm~\ref{alg3}. It is worth recalling that NbN training and PINN are unsupervised learning techniques for neural networks.
\begin{algorithm}[h!]
	\caption{Neuron-by-neuron training of neural networks with activation function $\rsigm$}\label{alg3}
	\KwData{---}
	\KwResult{Trained for the $E$ epoch NNs ${\vec u}_{\bm \theta}$, each of which consists of $N$ neurons on hidden layer}
	$\Delta x \gets X / N$\;
	${\Delta\zeta} \gets 1$\;
	Initialization of ${u}_{{\bm \theta};l}(x)$ of PINN ${\vec u}_{\bm \theta}$ with following steps\;
	\For{$k=0,\dots, N-1$}{
		${W}^{(1)}_k \gets 2 \Delta \zeta / \Delta x$\;
		${b}^{(1)}_k \gets -2 k \Delta \zeta$\;
		${W}^{(2)}_k \gets - \Delta x \mathcal{N}_l [{\vec u}(0), x_{k}] / \Delta \zeta$\;}
	${b}^{(2)}_0 \gets {u_l}(0)$\;
	Neuron-by-neuron improving weights ${W}^{(2)}_{k}$ with following steps\;
	\For{$e=0,\dots, E-1$}{
		\For{$k=0,\dots, N-1$}{
			${W}^{(2)}_{k} =  - \Delta x \mathcal{N}_l [{\vec u}_{\bm \theta}(x_{k}), x_{k}] / \Delta \zeta$.}
	}
	In our experiments, it was enough to use 3 epochs ($E=3$) of training to obtain high-precision results.
\end{algorithm}

\section{Numerical experiments}
For our experiments we used Pytorch~\cite{Paszke2019} version 2.1.2 and the training was carried out on a node with a GPU with characteristics similar to industry-leading GPU and CPU.


In the examples given in this section, we compare the accuracy of solutions predicted by a neural network using a sigmoid activation function (hereinafter referred to as $\sigma$) and a neural network that uses an activation function $\rsigm$. We use Algorithms~\ref{alg1} and ~\ref{alg2} for PIDD initialization of neural networks with $\sigma$ and $\rsigm$ activation functions, respectively. For all experiments, we use a uniform distribution for the coordinates of the collocation points. The reference solutions were obtained using the \verb*|odeint| solver of \verb*|scipy.integrate| library. As a measure of accuracy, we use the relative ${\mathbb L}_2$ error, which is defined in Appendix~\ref{appC}.

\subsection{Set of physical problems}
To compare the accuracy of neural networks, we consider the physical problems given in this subsection which are described by linear, nonlinear, and chaotic systems of differential equations.

\subsubsection{Harmonic Oscillator}
The harmonic oscillator is governed by a system of two hidden-order equations for the $t\in [0, T]$
\begin{eqnarray}
	&& \frac{d u_1}{d t} = u_2,\quad \frac{d u_2}{d t} = - \omega^2 u_1, \quad t \in [0,T],\label{eq58}\\
	&& u_1(0) = 1,\quad u_2(0) = 0 \label{eq59},
\end{eqnarray}
where $\omega$ is the frequency of the considered system. We used the following parameters: $T=100$ and $\omega = 1$. The exact analytical solution of this problem is $u_1^{\text{exact}} = \cos\left(\omega t\right)$, $u_2^{\text{exact}} = - \sin\left(\omega t\right)$.

\subsubsection{Relativistic slingshot}
Consider the problem of a source for single circularly polarized attosecond x-ray pulses~\cite{Wang2013}. The problem is formulated by following a system of ordinary differential equations in the absolute coordinate $t$ (see details in~\cite{eskin2024twohidden})
\begin{align}
	& \frac{\partial h}{\partial t} = \left(E_x - \varepsilon \frac{u_{\perp}^{2}}{1+u_{\perp}^{2}}\right) \frac{1}{1 + b},   \notag\\
	&\frac{\partial x}{\partial t} = \frac{b}{1 + b},\notag\\
	& \frac{\partial y}{\partial t} = \frac{u_y}{h (1+b)},\notag\\
	& \frac{\partial z}{\partial t} = \frac{u_z}{h (1+b)},\notag\\
	& b = \frac{1 + u_{\perp}^{2} - h^{2}}{2 h^{2}},\notag\\
	& u_{y} = a_{y,L} - \varepsilon y,\notag\\
	& u_{z}=a_{z,L}-\varepsilon z,\label{eq17}		
\end{align}
where $h=\gamma - u_{x}$, $\gamma^{2} = 1 + u_{\perp}^{2} + u_{x}^{2}$ is the relativistic gamma factor, $u_{x,y,z}$ are the space components of the four-velocity, $u_{\perp}^{2} = u_{y}^{2}+ u_{z}^{2}$, $\varepsilon = 2 \pi e^{2}n'l'/m\omega_{L}c$ ($l'$ is foil thickness), $a_{y,L}$ and $a_{z,L}$ are the $y$ and $z$ components of the laser vector potential ${\vec A}$. It is assumed that the electron is initially at $x=0$, $y=0$ and $z=0$, $E_{x} = \varepsilon {\rm th}\left({x}/ {4l'}\right)$. For the calculations we take $\varepsilon=4\pi$,  $l'=0.01 \lambda_{L}$ ($\lambda_{L} = 2 \pi {c}/{\omega_{L}}$).The period of impulse is $T_L = 1$ ($\omega_{L} = 2\pi / T_{L}$) and duration of impulse is  $T = 4 T_{L}$. In our calculation the value $l'$ was given $l'=0.01 2 \pi$. The amplitude of the vector potential components of the initiating pulse is described by
\begin{align}
	& a_{y} = a_0 \sin\left( 2 \pi t / T_L \right) \sin^2 \left( \pi t / T \right),\quad t \in [0,T],   \notag\\
	& a_{z} = a_{0} \cos\left( 2 \pi {t}/{T_{L}} \right) \sin^2 \left( \pi t / T \right), \quad t\in [0,T].
\end{align}

\subsubsection{Lorentz system}
As a further example, let us consider the chaotic Lorenz system. This system of equations arises in studies of convection and instability in planetary atmospheric convection, in which variables describe convective intensity and horizontal and vertical temperature differences~\cite{Lorenz1963DeterministicNF}. This system is given by the following set of ordinary differential equations:
\begin{align}
	&\frac{dx}{d t} = \tilde{\sigma} \left(y - x\right),\notag\\
	& \frac{dy}{d t} = x \left(\rho - z\right) - y,\notag\\
	& \frac{dz}{d t} = x y - \beta z,\label{eq54}		
\end{align}
where $\rho$, $\tilde{\sigma}$, and $\beta$ are the Prandtl number, Rayleigh number, and a geometric factor, respectively. We take the classical parameters $\tilde{\sigma} = 10$, $\rho = 28$ and $\beta = 8/3$. The max time is $T = 20$, and initial conditions are $x(0)=1$, $y(0)=1$, and $z(0)= 1$.

\subsection{Results of numerical experiments for the PIDD initialization}

In our experiments number of neurons of the hidden layer of neural networks is 20000, as the number of collocation points equals 20000. Results of experiments are presented in the Table~{\ref{table1}}. Column 3 (relative ${\mathbb L}_2$ errors) shows the relative errors between solutions, provided by neural networks and solutions obtained using \verb*|odeint| solver of \verb*|scipy.integrate| library ($u_1$ and $u_2$ for harmonic oscillator; $h$, $x$, $y$ and $z$ for relativistic slingshot; $x$, $y$ and $z$ for Lorentz system). The data presented herein provides compelling evidence that substituting the sigmoid activation function with the rectified sigmoid function results in a reduction of the relative error by at least an order of magnitude for any given value. Moreover, the duration of PIDD initialization is taken short time. 

Note, that calculating the relative errors between the exact solution of the harmonic oscillator problem and the solution by \verb*|odeint| for a selected set of points, we have the following values: $\epsilon[{u}^{\text{exact}}_1, {u^{\text{ref}}_1}] = 1.03{\times} 10^{-6}$, $\epsilon[{u}^{\text{exact}}_2, {u^{\text{ref}}_2}] = 1.04{\times} 10^{-6}$. Thus, these values are very close to the values obtained by a neural network with a rectified sigmoid activation function.

\begin{table*}[h!]
	\centering
	\begin{tabular}{cc|cc}
		\toprule
		\multicolumn{2}{c}{\bf{Experiments}} & \multicolumn{2}{c}{\bf{Performance}}\\
		\midrule
		\bf{Problem} & \bf{Activation}  &  \bf{Relative ${\mathbb L}_2$ errors} & \bf{Run time (s)} \\
		& \bf{function}  &  \bf{ ($\epsilon[{u}_{\bm \theta;l}, {u^{\text{ref}}_l}]$)} & \\
		\midrule
		{\bf Harmonic Oscillator} & $\sigma$ & {($5.67{\times} 10^{-5}$, $6.82{\times} 10^{-4}$)} & 0.0020 \\
		{\bf Harmonic Oscillator} & $\rsigm$ & {(${3.88{\times} 10^{-6}}$, ${2.63{\times} 10^{-6}}$)} & 0.0017 \\
		{\bf Relativistic slingshot} & $\sigma$ & {(${1.88{\times} 10^{-5}}$, ${1.92{\times} 10^{-5}}$, ${3.48{\times} 10^{-5}}$, ${3.43{\times} 10^{-5}}$)} & 0.004 \\
		{\bf Relativistic slingshot} & $\rsigm$ & {(${4.80{\times} 10^{-7}}$, ${3.57{\times} 10^{-7}}$, ${2.52{\times} 10^{-6}}$, ${2.70{\times} 10^{-6}}$)} & 0.004 \\
		{\bf Lorentz system} & $\sigma$ & {(${5.52{\times} 10^{-5}}$, ${3.07{\times} 10^{-4}}$, ${4.27{\times} 10^{-5}}$)} & 0.0026 \\
		{\bf Lorentz system} & $\rsigm$ & {(${5.14{\times} 10^{-6}}$, ${8.10{\times} 10^{-6}}$, ${3.03{\times} 10^{-6}}$)} & 0.0025 \\
		\bottomrule
	\end{tabular}
	\caption{Relative ${\mathbb L}_2$ error and run time of the PIDD initialization for the different physical problems and both neural networks with sigmoid activation function and $\rsigm$ activation function.}
	\label{table1}
\end{table*}

\subsection{Results of numerical experiments for the NbN training}
The whole time domains $[0;T]$ for all problems were split into 20 disjoint equivalent time windows of size $\Delta t = T / 20$. In our experiments for each time window we taken neural networks with 10000 neurons on the hidden layer. Results of experiments are presented in Table~{\ref{table2}} (the values are the same as in Table~\ref{table1}). It is evident from this table that neural networks with a rectified sigmoid activation function, trained using the NbN method exhibit superior accuracy compared to neural networks utilizing a sigmoid activation function. As NbN represents an unsupervised approach to training, it requires more time compared to PIDD initialization but less than vanilla PINN techniques~\cite{Raissi2019,Wang2022,eskin2024twohidden}.

\begin{figure}[t!]
	\centering
	\begin{subfigure}[t]{0.24\textwidth}
		\centering
		\includegraphics[width=\linewidth]{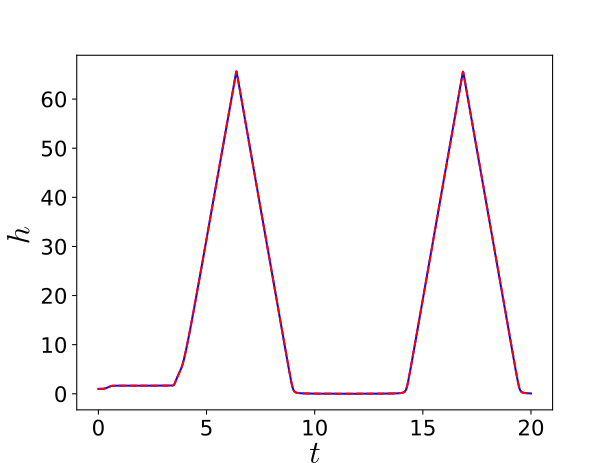}
		\caption{$h$}
	\end{subfigure}
	\begin{subfigure}[t]{0.24\textwidth}
		\centering
		\includegraphics[width=\linewidth]{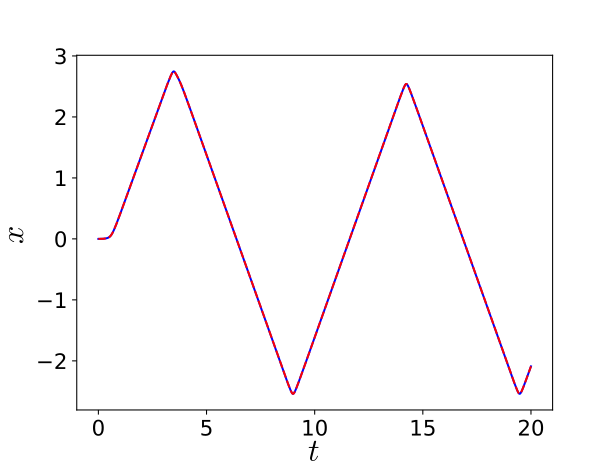}
		\caption{$x$}
	\end{subfigure}
	\begin{subfigure}[t]{0.24\textwidth}
		\centering
		\includegraphics[width=\linewidth]{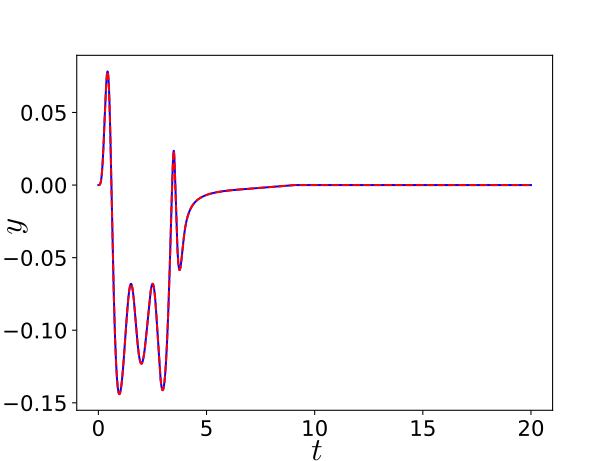}
		\caption{$y$}
	\end{subfigure}
	\begin{subfigure}[t]{0.24\textwidth}
		\centering
		\includegraphics[width=\linewidth]{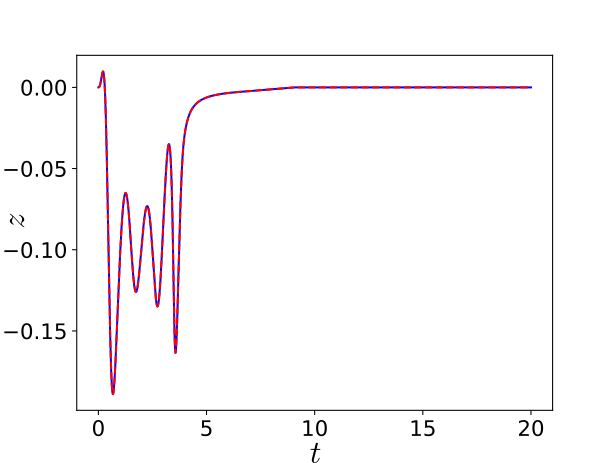}
		\caption{$z$}
	\end{subfigure}

	\begin{subfigure}[t]{0.24\textwidth}
		\centering
		\includegraphics[width=\linewidth]{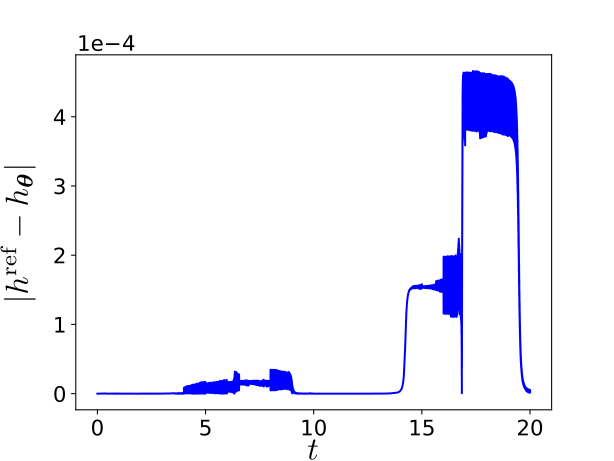}
		\caption{$h$}
	\end{subfigure}
	\begin{subfigure}[t]{0.24\textwidth}
		\centering
		\includegraphics[width=\linewidth]{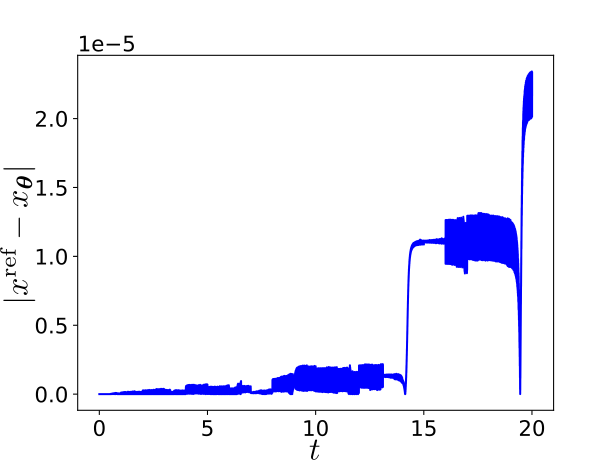}
		\caption{$x$}
	\end{subfigure}
	\begin{subfigure}[t]{0.24\textwidth}
		\centering
		\includegraphics[width=\linewidth]{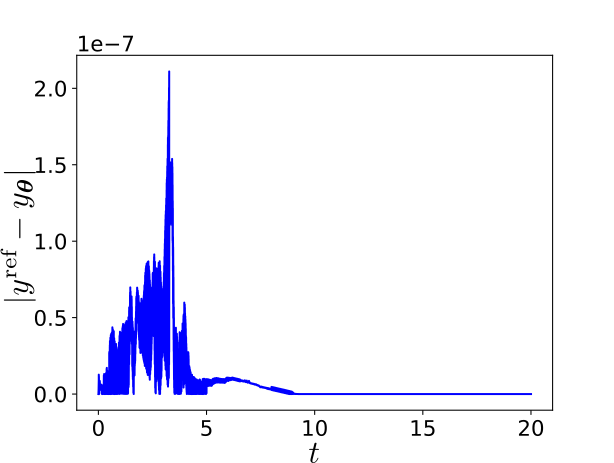}
		\caption{$y$}
	\end{subfigure}
	\begin{subfigure}[t]{0.24\textwidth}
		\centering
		\includegraphics[width=\linewidth]{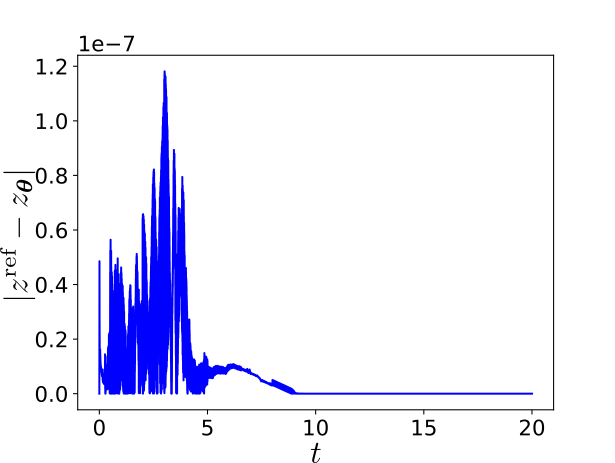}
		\caption{$z$}
	\end{subfigure}
	\caption{(a), (b), (c) and (d) are comparison of the predicted (red dash lines) and reference solutions (blue solid lines) corresponding to $h$, $x$, $y$ and $z$, respectively. (e), (f), (g) and (h) are absolute errors $|h^{\text{ref}} - h_{\bm \theta}|$, $|x^{\text{ref}} - x_{\bm \theta}|$, $|y^{\text{ref}} - y_{\bm \theta}|$, and $|z^{\text{ref}} - z_{\bm \theta}|$, respectively. }
	\label{fig1}
\end{figure}

Results of NbN training for neural networks with rectified sigmoid function as an activation function are shown on the Figure~\ref{fig1}. Over the all-time interval, the absolute errors for all values obtained by such neural networks by NbN training are significantly lower than the absolute errors obtained by the vanilla PINN~\cite{eskin2024twohidden}.

\begin{table*}[h!]
	\centering
	\begin{tabular}{cc|cc}
		\toprule
		\multicolumn{2}{c}{\bf{Experiments}} & \multicolumn{2}{c}{\bf{Performance}}\\
		\midrule
		\bf{Problem} & \bf{Activation}  &  \bf{Relative ${\mathbb L}_2$ errors} & \bf{Run time (s)}\\
		& \bf{function}  &  \bf{ ($\epsilon[{u}_{\bm \theta;l}, {u^{\text{ref}}_l}]$)} & \\
		\midrule
		{\bf Harmonic Oscillator} & $\sigma$ & {(${3.35\times 10^{-5}}$, ${3.17\times 10^{-5}}$)} & 302 \\
		{\bf Harmonic Oscillator} & $\rsigm$ & {(${3.35\times 10^{-6}}$, ${3.30\times 10^{-6}}$)} & 254 \\
		{\bf Relativistic slingshot} & $\sigma$ & {(${6.40\times 10^{-5}}$, ${5.19\times 10^{-5}}$, ${5.44\times 10^{-5}}$, ${6.08\times 10^{-5}}$)} & 982 \\
		{\bf Relativistic slingshot} & $\rsigm$ & {(${5.79\times 10^{-6}}$, ${4.18\times 10^{-6}}$, ${5.16\times 10^{-7}}$, ${3.16\times 10^{-7}}$)} & 903 \\
		{\bf Lorentz system} & $\sigma$ & {(${3.47\times 10^{-3}}$, ${5.04\times 10^{-3}}$, ${2.12\times 10^{-3}}$)} & 1077 \\
		{\bf Lorentz system} & $\rsigm$ & {(${2.63\times 10^{-3}}$, ${3.83\times 10^{-3}}$, ${1.62\times 10^{-3}}$)} & 888 \\
		\bottomrule
	\end{tabular}
	\caption{Relative ${\mathbb L}_2$ error and run time of the NbN training for the different physical problems and both neural networks with sigmoid activation function and $\rsigm$ activation function.}
	\label{table2}
\end{table*}

\section{Conclusions}

In this paper, a study has been carried out on neural networks with one hidden layer and a rectified sigmoid activation function for solving physical problems. Algorithms of physics-informed data-driven (PIDD) initialization of a neural network and gradient-free neuron-by-neuron (NbN) fitting method for neural networks with one hidden layer and rectified sigmoid function have been proposed. The results show that the accuracy of the ODE solutions which is provided by neural networks with a rectified sigmoid activation function outperforms the  accuracy of the ODE solution which is provided by neural networks with a sigmoid activation function.  This superiority is observed for the neural networks both after PIDD initialization and NbN training, for linear, nonlinear, and chaotic problems. Using the example of a harmonic oscillator, it is shown that the accuracy of solving the ODE using a neural network with such an activation function as a rectified sigmoid is close in values to standard numerical solvers. Our research continues to unlock the potential of neural networks with one or two hidden layers, and we hope that it will allow us to develop more accurate, effective, and reliable methods to solve complex physical problems for them.

\newpage

\appendix
\section{Physics-informed data-driven (PIDD) initialization of neural networks}\label{appA}
\begin{algorithm}[th!]
	\caption{Physics-informed data-driven initialization of neural networks}\label{alg1}
	\KwData{$\left\{{\vec u}_k\right\}_{k=0}^N$ for $\left\{x_k\right\}_{k=0}^N$ (uniform grid with a step $\Delta x \gets X / N$)}
	\KwResult{Initialized neural network ${u}_{{\bm \theta};l}(x)$ of PINN ${\vec u}_{\bm \theta}$, which consists of $N$ neurons on hidden layer}
	${\Delta\zeta} \gets \ln\left(2 + \sqrt{3}\right) / 2$\;
	\For{$k=0,\dots, N-1$}{
		${W}^{(1)}_k \gets 2 \Delta \zeta / \Delta x$\;
		${b}^{(1)}_k \gets -2 k \Delta \zeta$\;
		$\kappa_k \gets \sum\limits_{m = k-L}^{k+L} \sigma'\left({2\Delta \zeta} \left[k - m)\right]\right)$\;
		${W}^{(2)}_k \gets - \dfrac{\Delta x}{2\Delta \zeta} \dfrac{\mathcal{N}[{\vec u}_{k}, x_{k}]}{\kappa_k}$\;}
	${b}^{(2)}_0 \gets {u}_{l;0} - \sum\limits_{k = 0}^{N-1} {W}^{(2){\rm  T}}_k\sigma\left({W}^{(1){\rm T}}_k x_0 + {b}^{(1)}_k\right)$.
\end{algorithm}

\section{Neuron-by-neuron training of neural networks}\label{appB}

\begin{algorithm}[th!]
	\caption{Neuron-by-neuron training}\label{alg5}
	\KwData{---}
	\KwResult{Trained for the $E$ epoch NNs ${u}_{{\bm \theta};l}(x)$ of PINN ${\vec u}_{\bm \theta}$, each of which consists of $N$ neurons on hidden layer}
	$\Delta x \gets X / N$\;
	${\Delta\zeta} \gets \ln\left(2 + \sqrt{3}\right) / 2$\;
	Initialization of ${u}_{{\bm \theta};l}(x)$ with following steps\;
	\For{$k=0,\dots, N-1$}{
		${W}^{(1)}_k \gets 2 \Delta \zeta / \Delta x$\;
		${b}^{(1)}_k \gets -2 k \Delta \zeta$\;
		${W}^{(2)}_k \gets - 2 \Delta x \mathcal{N}_l [{\vec u}(0), x_{k}] / \Delta \zeta$\;}
	${b}^{(2)}_0 \gets {u_l}(0)$\;
	Neuron-by-neuron improving weights ${W}^{(2)}_{k}$ with following steps\;
	\For{$e=0,\dots, E-1$}{
		\For{$k=0,\dots, N-1$}{
			${W}^{(2)}_{k} =  - 2 \Delta x \mathcal{N}_l [{\vec u}_{\bm \theta}(x_{k}), x_{k}] / \Delta \zeta$.}
	}
	In our experiments, it was enough to use 3 epochs ($E=3$) of training to obtain high-precision results.
\end{algorithm}

\section{Measure of accuracy}\label{appC}
To evaluate the accuracy of the approximate solution obtained with the help of PINN method, the values of the solution of (\ref{eq1}) predicted by the neural network at given points are compared with the values calculated on the basis of classical high-precision numerical methods. As a measure of accuracy, the relative total ${\mathbb L}_2$ error of prediction is taken, which can be expressed with the following relation
\begin{equation}\label{eq61}
	\epsilon[{u}_{\bm \theta}, {u}] = \left\{\frac{1}{N_e} \sum^{N_e}_{i=1} \left[{u}_{\bm \theta}({\vec x}_i) - {u}({\vec x}_i)\right]^2 \right\}^{1/2} {\times} \left\{\frac{1}{N_e} \sum^{N_e}_{i=1} \left[{u}({\vec x}_i)\right]^2 \right\}^{-1/2},
\end{equation}
where $\left\{{\vec x}_{i}\right\}^{N_{e}}_{i=1}$ is the set of evaluation points taken from the domain $\Omega$, ${u}_{\bm \theta}$ and ${u}$ are the predicted and reference solutions respectively.

\newpage
\bibliographystyle{IEEEtran}
\bibliography{Eskin_rSigmoid}

\end{document}